\documentclass[12pt]{article}
\usepackage{a4,amssymb,amsmath}

\newtheorem{lem}{Lemma}
\newcommand{\bx}{\hspace*{\fill}$\square$}
\begin{document}
\noindent
{\Large \bf Pluripolar graphs are holomorphic}\\
\\
{\large \bf Nikolay Shcherbina}\\
\\
{\small Department of Mathematics, University of Wuppertal, 
D-42097  Wuppertal, Germany}\\
({\small e-mail: \tt{shcherbina@math.uni-wuppertal.de}})\\

\paragraph{{\large Abstract.}} Let $\Omega$ be a domain in ${\mathbb C}^n$ and let
$f: \Omega \rightarrow {\mathbb C}$ be a continuous function. We prove
that the graph $\Gamma(f)$ of the function $f$ is a pluripolar subset
of ${\mathbb C}^{n+1}$ if and only if  $f$ is holomorphic.

\section{Introduction}
A function $\varphi$ defined on a domain $U \subset
{\mathbb C}^n$ with values in $[-\infty,+\infty)$ is called {\em
  plurisubharmonic} in $U$ if $\varphi$ is
upper semicontinuous and its restriction to the components of the
intersection of a complex
line with $U$ is subharmonic.

A set $E \subset {\mathbb C}^n$ is called {\em pluripolar} if there is
a neighbourhood $U$ of $E$ and a plurisubharmonic function $\varphi$
on $U$ such that $E \subset \{\varphi=-\infty\}$. By a result of
B. Josefson [J], the function $\varphi$ in this definition can be
chosen to be plurisubharmonic in the whole of $\,{\mathbb C}^n$
(i.e. $U={\mathbb C}^n$).

In 1963 T. Nishino raised the following question in connection to his
paper [N1]:

\vskip0.7cm

{\em Let $\Delta$ be the unit disc in ${\mathbb C}_z$ and let $f: \Delta
  \rightarrow {\mathbb C}_w$ be a continuous function such that its
  graph $\Gamma (f)$ is a pluripolar subset of $\,{\mathbb
  C}^2_{z,w}$. Does it follow that $f$ is holomorphic?}

\vskip0.7cm

The main result of this paper gives a positive answer to Nishino's
question and can be formulated as follows.

\paragraph{Theorem.} {\em Let $\,\Omega$ be a domain in ${\mathbb C}^n$
  and let $f: \Omega \rightarrow {\mathbb C}$ be a continuous
  function. The graph $\Gamma (f)$ of the function $f$ is a pluripolar
  subset of $\,{\mathbb C}^{n+1}$ if and only if $f$ is
  holomorphic.}

\vskip0.7cm

As a consequence of Theorem one can easily obtain the following more 
general statement.

\paragraph{Corollary.} {\em Let $\,\Omega$ be a domain in ${\mathbb C}^n_z$
  and let $E$ be a closed subset of 
$\,\Omega\times {\mathbb C}_w \subset {\mathbb C}^{n+1}_{z,w}$ 
such that the fibers $E(z)=\{w\in \mathbb C_w:(z,w)\in E\}$ of
$E$ are finite and depend continuously on $z \in \Omega$ in the
Hausdorff metric. Assume that the number $\,\#E(z)$ of points in the
fiber $E(z)$ is bounded from above in $\Omega$. Then $E$ is  a pluripolar subset of 
$\,{\mathbb C}^{n+1}_{z,w}$ if and only if it has the form
\begin{equation}\label{1}
E=\{(z,w)\in \Omega\times {\mathbb C}_w: w^m+a_1(z)w^{m-1}+...+a_m(z)=0\}, 
\end{equation}
where the functions $a_1(z), a_2(z),..., a_m(z)$ are holomorphic in $\Omega$.}

\vskip0.7cm

Note that the proof of Theorem can not be directly applied to the set
$E$ described in Corollary. Namely, the topological argument used in the
proof of Lemma 3 and based on the fact that the first homology group 
$H_1(\Omega\times {\mathbb C}_w\setminus {\Gamma (f)}, {\mathbb  Z})$
is nontrivial does not work in this case. In the last 
section of the paper we construct an example of a compact subset 
$E$ of $\,\bar{\Delta}\times {\mathbb C}_w \subset {\mathbb C}^{2}_{z,w}$
 $(\Delta=\{z: |z|<1\})$ with finite fibers $E(z)$ depending 
continuously on $z\in \bar{\Delta}$ in the Hausdorff metric such that 
$H_1(\Delta\times {\mathbb C}_w\setminus E, {\mathbb  Z})=0$. 
In particular, there is a neighbourhood $U(E)$ of $E$ which 
does not contain any subset of $\bar{\Delta}\times {\mathbb C}_w$ defined by a 
Weierstrass pseudopolynomial (i.e. defined by the equation (1) with
$a_1(z), a_2(z),..., a_m(z)$ being continuous functions in $\Omega$).

\paragraph{Remark.}  In the special case when the function $f$ is assumed to be
$C^{1}$-smooth and its graph $\Gamma (f)$ is assumed to be complete
pluripolar (i.e. $\Gamma (f)=\{\varphi=-\infty\}$ for some function
$\varphi$, plurisubharmonic in a neighbourhood of $\Gamma (f)$), 
 a positive answer to  Nishino's question was given by Ohsawa [O] using
$L^{2}$ estimates for $\bar{\partial}$. In this case one can also
apply Pinchuk's method adapted to $C^1$-surfaces in [\v CH, p. 59-62] and
construct, to get a contradiction, a
one-parameter family of holomorphic disks $\{D_{\alpha}\}$ attached to
a totally real piece of $\Gamma (f)$ by an arc on the
boundary. Restricting the plurisubharmonic function $\varphi$ such
that $\Gamma (f) \subset \{\varphi=-\infty\}$ to each of these disks,
we get $\varphi \equiv  -\infty$ on $D_{\alpha}$ and, hence,
$\textstyle \bigcup_{\alpha}{D_{\alpha}} \subset \{\varphi=-\infty\}$
which gives the desired contradiction, since the set  $\textstyle
\bigcup_{\alpha}{D_{\alpha}}$ has real dimension $3$. Note that
neither of the methods
mentioned here can be applied to prove our Theorem.

\vskip0.9cm

{\small {\em Acknowledgements.} Part of this work was done while the
  author was a visitor at the Max Planck Institute of Mathematics
  (Bonn). It is my pleasure to thank this institution for its
  hospitality and excellent working conditions. I would like to thank
  E.M. Chirka who communicated to me the problem stated above,
  T. Ohsawa for informing me that the problem was first raised in 1963
  by T.Nishino, and E.L.Stout for pointing out to me the reference for 
the paper [A].}

\section{Preliminaries}

For bounded nonempty sets $E_1$ and $E_2$ in ${\mathbb C}_w$, the 
{\em Hausdorff distance} is defined as 
\[
d(E_1, E_2)={\sup_{w_2 \in E_2}\inf_{w_1 \in E_1}|w_1-w_2|}+{\sup_{w_2 \in E_1}\inf_{w_1 \in E_2}|w_1-w_2|}.
\]
A family of compact sets $E(z)$ in ${\mathbb C}_w$ parametrized by $z
\in \Omega \subset {\mathbb C}^n_z$ is said to be continuously
dependent on $z$ in the Hausdorff metric if, for each sequence
$\{z_n\}$ of points in $\Omega$ converging to a point $z_0 \in \Omega$,
one has $d(E(z_n), E(z_0)) \rightarrow 0$ as $n \rightarrow
\infty$. In particular, if $\Omega$ is a domain in ${\mathbb C}^n_z$
and $E$ is a nonempty closed subset of $\,\Omega\times {\mathbb C}_w$
with bounded fibers $E(z)=\{w\in {\mathbb C}_w : (z,w)\in E\}$
depending continuously on $z \in \Omega$ in the Hausdorff metric, then
each fiber $E(z)$, $z \in \Omega\,$, is nonempty.

\vskip0.3cm

For a compact set $K$ in ${\mathbb C}^n$, the {\em polynomial hull}
$\hat{K}$ of $K$ is defined as
\[
\hat{K}=\{z\in {\mathbb C}^n: |P(z)|\leq \sup_{w \in K} |P(w)|\mbox{ 
  for all holomorphic polynomials $P$ in ${\mathbb C}^n$}\}.
\]
The set $K$ is called {\em polynomially convex} if $\hat{K}=K$.

\vskip0.3cm

The first simple lemma is classical and follows, for example, from the
theorem 4.3.4 in [H].

\begin{lem}\label{lem1}
A compact set $K$ in ${\mathbb C}^n$ is polynomially convex if and
only if for any point $Q \in {\mathbb C}^n \setminus K$ there is a function $\varphi$,
plurisubharmonic in ${\mathbb C}^n$, such that
\begin{equation}\label{1}
\sup_{z\in K}\varphi(z)<\varphi(Q).
\end{equation}
\end{lem}

\begin{lem}\label{lem2}
Let $K$ be a polynomially convex compact set in ${\mathbb C}^n$ and let
$E$ be a pluripolar compact set in ${\mathbb C}^n$. Then the set
$(\widehat{K\cup E})\setminus K$ is pluripolar.
\end{lem}

\paragraph{Proof.} From pluripolarity of the set $E$ it follows that
there is a function $\varphi_E$,
plurisubharmonic in ${\mathbb C}^n$, such that $E \subset \{\varphi_E
= -\infty\}$. To prove Lemma 2, we
shall prove that $(\widehat{K\cup E})\setminus K \subset
\{\varphi_E=-\infty\}$.

Assume, by contradiction, that there is a point $Q\in
(\widehat{K\cup E})\setminus K$ such that $\varphi_E(Q)>-\infty$. Since $Q
\notin K$, and since the set $K$ is polynomially convex, it follows
from Lemma 1 that there is a
function $\varphi_K$, plurisubharmonic in ${\mathbb C}^n$, such that $\sup_{z \in K}
\varphi_K(z)<\varphi_K(Q)$. Then, for $\varepsilon$ positive and small
enough, one also has that $\sup_{z \in K}(\varphi_K(z)+\varepsilon
\varphi_E(z))<\varphi_K(Q)+\varepsilon \varphi_E(Q)$. Since
$\varphi_E (z)=-\infty$ for $z \in E$, it follows that $\sup_{z \in
  (K\cup E)}(\varphi_K(z)+\varepsilon
\varphi_E(z))<\varphi_K(Q)+\varepsilon \varphi_E(Q)$. By
Lemma 1 applied to the function $\varphi_K+\varepsilon \varphi_E$, we
get that $Q \notin (\widehat{K\cup E})$. This gives the desired
  contradiction. \bx

\vskip0.5cm

The next statement was first proved by H. Alexander (see Corollary 1
in [A]). For the
convenience of reading we include here its proof.

\vskip0.5cm

\begin{lem}\label{lem3}
Let $U$ be a bounded domain in ${\mathbb C}_z\times {\mathbb R}_u
\subset {\mathbb C}^2_{z,w}(w=u+iv)$ and let $g: bU\rightarrow {\mathbb R}_v$
be a continuous function. Then $U
\subset \pi (\widehat{\Gamma (g)})$, where $\Gamma (g)$ is the graph of
$g$ and $\pi: {\mathbb
  C}^2_{z,w}\rightarrow {\mathbb C}_z \times {\mathbb R}_u$ is the
projection.
\end{lem}

\paragraph{Proof.} Consider an approximation of the domain $U$ by an
increasing sequence $\{U_n\}$ of domains with smooth
boundary. Further, consider a sequence of smooth functions $\{g_n\}$,
$g_n: bU_n \rightarrow {\mathbb R}_v$, which approximate the function
$g$, i.e., $\Gamma (g_n) \rightarrow \Gamma (g)$ in the Hausdorff
metric. Then it follows from the definition of polynomial hull that
$\limsup_{n\rightarrow \infty}\widehat{\Gamma (g_n)} \subset 
\widehat{\Gamma (g)}$, where convergence is understood to be in the Hausdorff
metric. Hence, it is enough to prove the statement of Lemma 3 in the
case where the domain $U$ has a smooth boundary and the function $g$ is
smooth.

Now we argue by contradiction and suppose that there is a point $Q \in
U \setminus \pi (\widehat {\Gamma(g)})$. Without loss of generality we may
assume that $Q$ is the origin $O$ in ${\mathbb C}_z \times {\mathbb
R}_u$. We know by Browder [B] that $\check{H}^2(\widehat{\Gamma(g)}, {\mathbb
C})=0$ (here $\check{H}^2(\widehat{\Gamma(g)}, {\mathbb C})$ is the second 
\v{C}ech cohomology group with complex coefficients). Then, by
Alexander duality (see, for example [S], p.296), we get
$H_1({\mathbb C}^2_{z,w}\setminus \widehat{\Gamma (g)}, {\mathbb
  C})=\check{H}^2(\widehat{\Gamma(g)}, {\mathbb C})=0$ 
(here $H_1({\mathbb C}^2_{z,w}\setminus \widehat{\Gamma (g)}, {\mathbb
  C})$ is
the first singular homology group with complex coefficients). On the
other hand, since $O \in U \setminus \widehat{\Gamma (g)}$, it follows
that the curve $\gamma_R$ consisting of the segment $\{z=0, u=0,
-R\leq v \leq R\}$ and the half-circle $\{z=0, w=\mbox{Re}^{i\theta},
-\frac{\pi}{2} \leq \theta \leq \frac{\pi}{2}\}$ do not intersect the
set $\widehat{\Gamma(g)}$ for $R$ big enough. Moreover, the linking
number of $\Gamma (g)$ and $\gamma_R$ is not equal to zero. Therefore,
$H_1({\mathbb C}^2_{z,w} \setminus \widehat{\Gamma (g)}, {\mathbb C})\neq
0$. This is a contradiction and the lemma follows. \bx

\vskip0.5cm

\begin{lem}\label{lem4}
Let $U$ be a simply connected domain in ${\mathbb C}_z$ and let
$f(z)=u(z)+iv(z): U\rightarrow {\mathbb C}_w$ be a function such that
both $u(z)$ and $v(z)$ are harmonic in $U$. If the graph $\Gamma
(f)$ of the function $f$ is a pluripolar subset of ${\mathbb
  C}^2_{z,w}$, then $f$ is holomorphic.
\end{lem}

\paragraph{Proof.} If $f$ is not holomorphic, we argue by
contradiction and suppose that
the set $\Gamma (f)$ is pluripolar. Then there is a  function
$\varphi$, plurisubharmonic
in ${\mathbb C}^2_{z,w}$, such that $\Gamma (f)
\subset \{\varphi=-\infty\}$. Let $\tilde{v}$ be the harmonic
conjugate function to $u$ in the domain $U$ such that
$\tilde{v}(z_0)=v(z_0)$ for some fixed point $z_0 \in U$. Then the
set $\{z \in U: \tilde{v}(z)+\varepsilon=v(z)\}$ is nonempty and
consists of real analytic curves for all $\varepsilon$ small
enough. Therefore, each of the holomorphic curves
$\Gamma_\varepsilon=\{(z,w): z\in U,
w=u(z)+i(\tilde{v}(z)+\varepsilon)\}$ intersects the set $\Gamma
(f)\subset \{\varphi=-\infty\}$ in real analytic curves. Since a real analytic
curve is not polar (see, for example [T, Th.II.26, p.50]), it follows that
$\Gamma_\varepsilon \subset \{\varphi=-\infty\}$ for all $\varepsilon$
small enough. This implies that $\varphi\equiv -\infty$ in ${\mathbb
  C}^2_{z,w}$ and gives the desired contradiction. \bx

\section{Proof of Theorem and Corollary} 

\paragraph{Proof of Theorem.} If the function $f$ is holomorphic, then 
the same argument as in the
proof of Lemma 4 shows that $\Gamma (f)$ is pluripolar. Namely, the
function
\[
\varphi(z_1,\ldots, z_{n+1})=\ln |z_{n+1}-f(z_1,\ldots, z_n)|
\]
is plurisubharmonic in $\Omega \times {\mathbb C}$ and $\Gamma
(f)=\{\varphi=-\infty\}$. Therefore, the set $\Gamma (f)$ is
pluripolar in ${\mathbb C}^{n+1}$.

Suppose now that the graph $\Gamma (f)$ of $f$ is pluripolar. To prove
that $f$ is holomorphic we consider two cases.

\paragraph{ {\it 1. The special case $n=1$.}} In this case $\Omega$ is a
  domain in ${\mathbb C}_z$ and $f(z)=u(z)+iv(z):\Omega \rightarrow
  {\mathbb C}_w$ is a continuous function such that its graph is
  pluripolar. Since holomorphicity is a local property, we can
  restrict ourselves to the case when $\Omega$ is a disc in ${\mathbb
  C}_z$ and, moreover, to simplify our notations, we can assume
  without loss of generality that $\Omega =\Delta=\{z: |z|<1\}$ is
  the unit disc and that the function $f$ is continuous on its closure
  $\bar{\Delta}$. It follows from Lemma 4 that either the function $f$
  is holomorphic or at least one of the functions $u$ and $v$ is not
  harmonic. Since both cases can be treated the same way, we can,
  to get a contradiction, assume that the function $u$ is not
  harmonic. Denote by $\tilde{u}$ the solution of the Dirichlet
  problem on $\Delta$ with boundary data $u$. Since $u$ is not
  harmonic, one has that $\tilde{u}\neq u$ in $\Delta$. Without
  loss of generality we can assume that
\begin{equation}\label{2}
u(z_0)<\tilde{u}(z_0)
\end{equation}
for some $z_0 \in \Delta$. Let
\[
C=\max \{\sup_{z\in \bar{\Delta}}|u(z)|, \sup_{z\in
  \bar{\Delta}}|v(z)|\}.
\]
Consider the set
\[
K=\{(z,w)\in \bar{\Delta}\times {\mathbb C}_w: \tilde{u}(z)\leq u \leq
3C, |v|\leq C\}.
\]

\vskip0.5cm

\begin{lem}\label{lem5}
The set $K$ is polynomially convex.
\end{lem}

\paragraph{Proof.} To prove polynomial convexity of $K$ we use Lemma
1. Consider an arbitrary point $(z^\ast, w^\ast) \in {\mathbb
  C}^2_{z,w}\setminus K$. If the point $(z^\ast, w^\ast)$ belongs to
  the set
\[
A_1=\{(z,w)\in {\mathbb C}^2_{z,w}:|z|>1 \mbox{ or } u>3C \mbox{ or }
|v|>C\},
\]
then inequality (2) will be satisfied for the point $Q=(z^\ast,
w^\ast)$ and the function
\[
\varphi_1(z,w)=\max\{|z|-1, u-3C, |v|-C\}
\]
plurisubharmonic in ${\mathbb C}^2_{z,w}$.

If the point $(z^\ast, w^\ast)$, $w^\ast=u^\ast+iv^\ast$ belongs to the set
\[
A_2=\{(z,w)\in \bar{\Delta}\times {\mathbb C}_w: u<\tilde{u}(z)\},
\]
then $u^\ast <\tilde{u}(z^\ast)$. Let $\varepsilon=\frac{1}{3}
(\tilde{u}(z^\ast)-u^\ast)$ and consider a function
$\tilde{u}_\varepsilon$ harmonic on the whole of ${\mathbb C}_z$ such
that $\max_{z \in
  \bar{\Delta}}|\tilde{u}(z)-\tilde{u}_\varepsilon(z)|<\varepsilon$. 
Since for $(z,w) \in K$ one has $u \geq \tilde{u}(z)\geq
\tilde{u}_\varepsilon(z)-\varepsilon$, and since 
$u^\ast=\tilde{u}(z^\ast)-3\varepsilon<\tilde{u}_\varepsilon(z^\ast)-2\varepsilon$, 
it follows that inequality (2) will be satisfied for the point
$Q=(z^\ast, w^\ast)$ and the function
\[
\varphi_2(z,w)=\tilde{u}_\varepsilon(z)-u
\]
plurisubharmonic in ${\mathbb C}^2_{z,w}$.

Since ${\mathbb C}^2_{z,w}\setminus K=A_1 \cup A_2$, we conclude from
Lemma 1 that the set $K$ is polynomially convex. This completes the
proof of Lemma 5. \bx

\vskip0.7cm

Consider now the domain
\[
U=\{(z,u) \in \Delta \times {\mathbb R}_u: u(z)<u<u(z)+2C\}
\]
in ${\mathbb C}_z \times {\mathbb R}_u$ and the real-valued function
$g(z,u)=v(z)$ on $bU$. Since $\sup_{z\in \bar{\Delta}}|u(z)|\leq C$,
one has $\sup_{z\in \bar{\Delta}}|\tilde{u}(z)|\leq C$ and hence
$\tilde{u}(z)\leq u(z)+2C\leq 3C$. It then follows from the
definitions of $U$ and $g$ that the graph $\Gamma (g)$ of the function
$g$ is contained in the set $\Gamma (f) \cup K$. Therefore, we get
$\widehat{\Gamma (g)} \subset (\Gamma \widehat{(f)\cup K})$. Since, by
Lemma 3, $\pi (\widehat{\Gamma (g)})\supset U$,  we conclude that
\begin{equation}\label{3}
\pi (\widehat{\Gamma (f) \cup K})\supset U.
\end{equation}

Consider the following open subset of $U$ :
\[
\tilde{U}=\{(z,u) \in \Delta \times {\mathbb R}_u:
u(z)<u<\tilde{u}(z)\}.
\]
Inequality (3) obviously implies that the set $\tilde{U}$ is
nonempty. Since, by the definition of the sets $K$ and $\tilde{U}$,
$\pi (K) \cap \tilde{U}=\emptyset$, it follows from (4) that
\begin{equation}\label{4}
\pi ((\widehat{\Gamma (f)\cup K})\setminus K)\supset \tilde{U}.
\end{equation}
Since, by our assumption, the graph $\Gamma (f)$ of $f$ is pluripolar,
we conclude from Lemma 2 and Lemma 5 that the set $(\widehat{\Gamma
  (f)\cup K})\setminus K$ is pluripolar, i.e.,
\begin{equation}\label{5}
(\widehat{\Gamma (f)\cup K}) \setminus K \subset \{\varphi=-\infty\}
\end{equation}
for some plurisubharmonic function $\varphi$.

From (3) one has that there is a neighbourhood $V$ of the point $z_0$
in ${\mathbb C}_z$ such that
\begin{equation}\label{6}
u(z)<\tilde{u}(z)
\end{equation}
for all $z \in V$. For each $a \in {\mathbb C}$ consider the complex
line $\ell_a=\{(z,w) \in {\mathbb C}^2: z=a\}$ and the set
\[
E_a=((\widehat{\Gamma (f)\cup K})\setminus K)\cap \ell_a.
\]
It follows from (5) and (7) that for $a \in V$ the projection of $E_a$
on the real line $\ell_a \cap \{v=0\}$ contains an open segment. Since
a polar set in ${\mathbb C}$ has Hausdorff dimension zero (see, for
example [T, Th.III.19, p.65]), it cannot be projected on an open segment in
${\mathbb R}$. Therefore, the set $E_a$ is not polar. It then
follows from (6) that $\varphi\equiv -\infty$ on $\ell_a$. Since this
argument holds true for all $a \in V$, we conclude that $\varphi
\equiv -\infty$ on ${\mathbb C}^2_{z,w}$. This contradiction proves
Theorem in the case $n=1$.

\paragraph{{\it 2. The general case.}} Let $k$ be one of the numbers
$1,2,\ldots, n$. For each ${\mathbf a}=(a_1, a_2,\ldots, a_n)\in
\Omega$ consider the function
\[
f^{\bf a}_k(z_k)=f(a_1,\ldots, a_{k-1}, z_k, a_{k+1},\ldots, a_n)
\]
defined on the domain
\[
\Omega^{\bf a}_k=\Omega \cap \{z_1=a_1,\ldots, z_{k-1}=a_{k-1},
z_{k+1}=a_{k+1},\ldots, z_n=a_n\} \subset {\mathbb C}_{z_k}.
\]
Since, by our assumptions, the set $\Gamma (f)$ is pluripolar, there
is a function $\varphi$, plurisubharmonic in ${\mathbb C}^{n+1}$, such
that $\Gamma (f) \subset \{\varphi=-\infty\}$. For all points
${\mathbf a}$ except for a pluripolar set in ${\mathbb C}^n$ one
obviously has that the function
\[
\varphi^{\bf a}_k (z_k, z_{n+1})=\varphi(a_1,\ldots, a_{k-1}, z_k,
a_{k+1},\ldots, a_n, z_{n+1})
\]
is not identically equal to $-\infty$ in ${\mathbb C}^2_{z_k,
  z_{n+1}}$. For all such points ${\mathbf a}$ we can use the argument
  from case 1 and conclude from continuity of the function $f^{\bf
  a}_k: \Omega^{\bf a}_k\rightarrow {\mathbb C}_{z_{n+1}}$ and from
  the inclusion $\Gamma (f^{\bf a}_k)\subset \{\varphi^{\bf
  a}_k=-\infty\}$ that the function $f^{\bf a}_k$ is
  holomorphic. Since the complement of a pluripolar set is everywhere
  dense, it follows from continuity of $f$ that the functions $f^{\bf
  a}_k$ are holomorphic for all ${\mathbf a} \in \Omega$. This
  argument holds true for any $k=1,2,\ldots, n$, so we
  conclude from the classical Hartogs theorem on separate analyticity
  that the function $f$ is holomorphic. The proof of Theorem is
  now completed. \bx

\paragraph{Proof of Corollary.} Since, by our assumption, the number 
$\,\#E(z)$ of points in the fiber of $E$ is bounded from above
  in $\Omega$, we can consider $k=\max_{z \in
  \Omega} \#E(z)$ and then the open subset ${\mathcal U}=\{z \in
  \Omega: \,\#E(z)=k\}$ of $\Omega$. Let $z_0$ be a point of $\mathcal U$ and
  let $h_i(z)$, $i=1, 2,...,k$, be the functions defining single-valued 
branches of $E(z)$ in a neighbourhood $U$ of $z_0$. Since, by our
  assumption, $E(z)$ depends continuously on $z \in
  \Omega$ in the Hausdorff metric, we conclude from Theorem that
  the functions $h_i(z)$ are holomorphic in $U$. Hence,
  $F(z)=\Pi_{i \neq j} (h_i(z)-h_j(z))$ is a well defined holomorphic
  function in $\mathcal U$ such that for each $z{'} \in b\,{\mathcal U} \cap
  \Omega$ one has $F(z) \rightarrow 0$ as $z \rightarrow z{'}$, $z \in
  {\mathcal U}$. Then the function
$$
\tilde{F}(z)=\left\{
\begin{array}{ll}
F(z),&  \text{ for } \quad  z \in {\mathcal U},\\
0,& \text{ for } \quad  z \in \Omega \setminus {\mathcal U}
\end{array}
\right.
$$
is continuous in $\Omega$ and holomorphic in ${\mathcal U}=\Omega \setminus
\{z: \tilde{F}(z)=0\}$. Therefore, by Rad\'o's theorem (see, e.g. [\v
C], p. 302), $\tilde{F}$ is holomorphic in $\Omega$. In
particular, the set $\{z \in \Omega: \tilde{F}(z)=0\}$ is an analytic
hypersurface.

Consider now the function $\,\Pi_{i=1}^k
(w-h_i(z))=w^k+a_1(z)w^{k-1}+\, ...\, +a_k(z)$. Since
$a_1(z), a_2(z), .., a_k(z)$ are symmetric functions of $h_1(z),
h_2(z), ..., h_k(z)$, they are  well defined and holomorphic in ${\mathcal U}$.
Moreover, since $E(z)$ depends continuously on $z \in
  \Omega$ in the Hausdorff metric, these functions are locally bounded
  near the set $\Omega \setminus  {\mathcal U}= \{z: \tilde{F}(z)=0\} $. It
  follows then from removability of
analytic singularities that the functions 
$a_1(z), a_2(z), ..., a_k(z)$ are holomorphic in the whole of $\,\Omega$. Since, by
our construction,
$E=\{(z,w)\in \Omega\times {\mathbb C}_w: w^m+a_1(z)w^{m-1}+\, ...\, +a_m(z)=0\}$,
 the corollary follows.
\bx

\paragraph{Remark.} The statement of Corollary was first proved in
[Sh] for sets represented by Weierstrass pseudopolynomials by a
different (and more complicated) method. It was later observed
independently by the author and by A. Edigarian [E] that the methods of
chapter 4 in [N2] give a simpler proof for these sets.

\section{Example}

We first prove the following simple lemma.

\begin{lem}\label{lem1}
Let $f$ and $g$ be holomorphic functions, defined in a neighbourhood
$U$ of a point $a \in {\mathbb C}_z$, such that $f(a)=g(a)$ and
$f^{'}(a) \neq g^{'}(a)$. Let $r$ be a positive number such that 
${\bar{\Delta}_{r}(a)=\{z \in {\mathbb C}_z : |z-a| \leq r\}} \subset
U$ and $f(z) \neq g(z)$ for $z \in {\bar{\Delta}_{r}(a) \setminus
  \{a\}}$. Then for all sufficiently small $\varepsilon>0$ the complex 
curve $\Sigma \subset {{\Delta}_{r}(a) \times {\mathbb C}_w}$ defined by the equation

\begin{equation}\label{8}
G(z,w)\stackrel{\rm def}{=}(w-f(z))(w-g(z))-\varepsilon=0
\end{equation}
is a branched covering over the disk ${\Delta}_{r}(a)$ with two branches and two branching points 
\begin{equation}\label{9}
b^{\pm}=a \pm \frac{2i}{f^{'}(a)-g^{'}(a)}{\sqrt{\varepsilon}}+O(\varepsilon).
\end{equation}
\end{lem}
\vskip0.5cm

\paragraph{Proof.} Equation (8) is quadratic with respect to $w$,
hence $\Sigma$ is a branched covering over ${\Delta}_{r}(a)$ with two
branches. A point $b$ is a branching point of $\Sigma$ if for some
$w_b$ one has $0=G^{'}_{w}(b,w_b)=2w_b-f(b)-g(b)$. Therefore, 
$w_{b}=\frac{1}{2}{(f(b)+g(b))}$ and then (8) implies that 
$-\frac{1}{4}{(f(b)-g(b))}^{2}-\varepsilon=0$, i.e.

\begin{equation}\label{18}
f(b)-g(b)=\pm2i{\sqrt{\varepsilon}}.
\end{equation}
Hence, in view of our choice of $r$, $b \rightarrow a$ as $\varepsilon
\rightarrow 0$. Then, using Taylor expansions of $f$ and $g$ at the
point $a$, we conclude from (10) and the assumption $f(a)=g(a)$ that 
$(f^{'}(a)-g^{'}(a))(b-a)+O(|b-a|^2)=\pm2i{\sqrt{\varepsilon}}$.
Finally, the assumption $f^{'}(a) \neq g^{'}(a)$ implies that
\[
{b-a}={\pm \frac{2i}{f^{'}(a)-g^{'}(a)}{\sqrt{\varepsilon}}+O(|b-a|^2)}={\pm \frac{2i}{f^{'}(a)-g^{'}(a)}{\sqrt{\varepsilon}}+O(\varepsilon)}.
\]
\bx

\vskip0.5cm

\paragraph{Construction of the set $E$.} Let $\rho$ be a smooth
real-valued function defined on the segment $[0,1]$ such that

$$
\rho(t)=\left\{
\begin{array}{lll}
1,&  \text{ for } \quad  0 \leq t \leq \frac{1}{3},\\
\text{decreasing},& \text{ for } \quad  \frac{1}{3} < t < \frac{2}{3},\\
0,& \text{ for } \quad  \frac{2}{3} \leq t \leq 1.
\end{array}
\right.
$$
Consider the set 
\[
E_1=\{(z,w) \in {\bar{\Delta} \times {\mathbb C}_w} : w^2=\rho(|z|)z\},
\]
where, as above, $\Delta=\{z \in {\mathbb C}_z : |z|<1\}$ is the unit
disk. This set has two branches over the disk
$\Delta_{\frac{2}{3}}(0)$ with one branching point at $z=0$. The 
branches are glued to each other along the circle ${\mathcal
  A}=\{|z|=\frac{2}{3}, w=0\}$ and become one branch $\{w=0\}$ for 
$\frac{2}{3} \leq |z| \leq 1$. Consider some points $A_1=(a_1 , 0)$ 
and $A_3=(a_3 , \sqrt{a_3})$ of $E_1$ and a point $A_2=(a_2 , C)$ 
with $a_1$, $a_2$, $a_3$ and $C$ real and positive such that 
$\frac{2}{3} <a_1 <1$, $0 < a_3 <\frac{1}{3}$ and $a_3< a_2 <a_1$. 
Further, consider the complex line ${\mathcal L}^{'}$ passing through 
the points $A_2$ and $A_1$ and the complex line ${\mathcal L}^{''}$ 
passing through the points $A_2$ and $A_3$. Let $a_1$, $a_2$, $a_3$ 
be already chosen and consider $C$ so big that the line ${\mathcal
  L}^{''}$ intersects $E_1$ in two points $A_3$ and 
$A^{'}_3=(a^{'}_3,-{\sqrt{a^{'}_3}})$, with $a^{'}_3$ real such that 
$0 < a^{'}_3 < a_3$, and the line ${\mathcal L}^{'}$ intersects $E_1$ 
only at the point $A_1$. The set $E$ will be constructed as a small 
deformation of the set 
${E_1} \cup (({\mathcal L}^{'} \cup {\mathcal L}^{''}) \cap
({\bar{\Delta} \times {\mathbb C}_w}))$ 
near the points $A_k$, $k=1, 2, 3$, that creates, as in Lemma 6, 
two branching points instead of each self-intersection point.

Let $r>0$ be so small that the disks
${\bar{\Delta}}_1={\bar{\Delta}_r}(a_1)$,
${\bar{\Delta}}_2={\bar{\Delta}_r}(a_2)$ and
${\bar{\Delta}}_3={\bar{\Delta}_r}(a_3)$ neither intersect each other
nor the circle $\{|z|=\frac{2}{3}\}$ and, moreover, do not contain the
point $a^{'}_3$. Denote by ${\mathcal E}_1$ the set 
$(E_1 \cup {\mathcal L}^{'}) \cap  ({\Delta}_1 \times {\mathbb C}_w)$, 
by ${\mathcal E}_2$ the set 
$({\mathcal L}^{'} \cup {\mathcal L}^{''}) \cap  ({\Delta}_2 \times
{\mathbb C}_w)$ and by ${\mathcal E}_3$ the connected component of the 
set $(E_1 \cup {\mathcal L}^{''}) \cap  ({\Delta}_3 \times {\mathbb
  C}_w)$ containing the point $A_3$. Then each of the sets 
${\mathcal E}_k$, $k=1, 2, 3$, is the union of the graphs of two 
holomorphic functions $f^{j}_k$, $j=1, 2$, having the same value and 
different derivatives, both of them real (which is easy to check by 
direct calculation) at the center of the respective disk
${\Delta}_k$. Therefore, we can apply Lemma 6 to each of these sets
and, if $\varepsilon$ is small enough, we will get branched coverings 
${\Sigma}_1$, ${\Sigma}_2$ and ${\Sigma}_3$ over the disks 
${\Delta}_1$, ${\Delta}_2$ and ${\Delta}_3$, respectively, with two 
branches and two branching points contained in the smaller disks 
${\Delta}^{'}_1={{\Delta}_{\frac{r}{3}}}(a_1)$,
${\Delta}^{'}_2={{\Delta}_{\frac{r}{3}}}(a_2)$ and
${\Delta}^{'}_3={{\Delta}_{\frac{r}{3}}}(a_3)$. Moreover, since 
for each $k=1, 2, 3$ the derivatives at the centers of the disks
$\Delta_k$ of the functions $f^{j}_k$, $j=1, 2$, are real, we conclude 
from (9) that one of the two branching points contained 
in ${\Delta}^{'}_k$ is contained in the half-disk 
$\{z \in {\Delta}^{'}_k : \operatorname{Im}z>0 \}$, while the other 
is contained in the half-disk 
$\{z \in {\Delta}^{'}_k : \operatorname{Im}z<0 \}$. Since both
branching points of each set $\Sigma_k$ are contained in the
respective disk $\Delta^{'}_k$, the set 
$\Sigma_k \cap ((\Delta_k \setminus \Delta^{'}_k) \times {\mathbb
  C}_w)$ will be the union of the graphs of two holomorphic functions 
$\tilde{f}^{j}_k$, $j=1, 2$, defined on $\Delta_k \setminus
\Delta^{'}_k$ and, moreover, if $\varepsilon$ is small enough, then 
each  function $\tilde{f}^{j}_k$ will be close enough to the 
corresponding function $f^{j}_k$. Define the functions
\[
{\tilde{\tilde{f}}^{j}_k}(z)=\rho\left(\frac{|z-a_k|}{r}\right){{{\tilde{f}}^{j}_k}(z)+\left(1-\rho\left(\frac{|z-a_k|}{r}\right)\right){{f^{j}_k}(z)}},
\]
for $z \in \Delta_k \setminus \Delta^{'}_k$, $k=1, 2, 3$, $j=1,
2$. Let $\tilde{\Sigma}_k$ be the union of the graphs of 
$\tilde{\tilde{f}}^{1}_k$ and $\tilde{\tilde{f}}^{2}_k$.  Now we 
can define the set $E$ as
$$\textstyle
E=\left(\left(E_1 \cup \left(\left({\mathcal L}^{'} \cup {\mathcal L}^{''}\right) \cap \left(\bar{\Delta} \times {\mathbb C}_w\right)\right)\right) \setminus \bigcup\limits^3_{k=1} {\mathcal E}_k\right) \cup \bigcup\limits^3_{k=1}\left(\tilde{\Sigma}_k \cup \left(\Sigma_k \cap \left(\bar{\Delta}^{'}_k \times {\mathbb C}_w\right)\right)\right).
$$
Define also the set $E^{\operatorname{reg}}$ as $E$ with the circle 
${\mathcal A}$, the point $A^{'}_3$ of the transversal
self-intersection of $E$ and all the branching points of $E$ being
removed. Then, by our construction, $E^{\operatorname{reg}}$ is a 
smooth connected $2$-dimensional surface transversal to the $w$-direction.

Note that each fiber $E(z)$ of the set $E$ has at most four points and that the fibers $E(z)$ depend continuously on $z \in \bar{\Delta}$ in the Hausdorff metric.

\paragraph{Claim 1.} $H_1(\Delta\times {\mathbb C}_w\setminus E, {\mathbb  Z})=0$. \\

\vspace{-0.6truecm}

\paragraph{Proof.} Let $a$ be a real positive number such that $a_3
\leq a < \frac{1}{3}$. Consider the point $A=(a,-\sqrt{a}) \in E$ and
a disk $\bar{\mathcal D}_s=\{(z,w): z=a, |w+\sqrt{a}| \leq s\}$ so
small that it intersects the set $E$ only at the point $A$. We first prove that the circle ${\mathcal C}_s=b{\mathcal D}_s$ is homological to zero in $\Delta\times {\mathbb C}_w\setminus E$.

Consider the curve $z(t)$ in ${\mathbb C}_z$ defined as

$$
z(t)=\left\{
\begin{array}{lll}
a(1-t)+(a_1+r)t,&  \text{ for } \quad  0 \leq t \leq 1,\\
a_1+re^{{\pi}i(t-1)},& \text{ for } \quad  1 < t \leq 2,\\
(a_1-r)(3-t)+(a_3+r)(t-2),& \text{ for } \quad  2 < t \leq 3,\\
a_3+re^{{\pi}i(t-3)},& \text{ for } \quad  3 < t \leq 4,\\
(a_3-r)(5-t)+{\frac{2}{3}}(t-4),& \text{ for } \quad  4 < t \leq 5.\\
\end{array}
\right.
$$
If $\pi_z: {\mathbb C}^2_{z,w} \rightarrow {\mathbb C}_z$ is the projection, then the curve $z(t)$ admits a uniquely defined lifting by ${\pi}^{-1}_z$ to the piece-wise smooth curve ${\gamma}$ in $E$ with the initial point $A$.

The curve $\gamma$ is transversal to the $w$-direction and has one
point of self-intersection, namely, the endpoint $(\frac{2}{3}, 0)$, where two smooth curves on the side $\{|z|<\frac{2}{3}\}$ meet each other. 

The geometric description of the curve $\gamma$ looks as follows. We start from the point $A=(a,-\sqrt{a})$ and then, over the segment $\{a \leq \operatorname{Re}z < \frac{2}{3}, \operatorname{Im}z=0\}$, the curve $\gamma$ is contained in  the ``lower'' branch of the set $E_1$, while over the segment $\{\frac{2}{3} \leq \operatorname{Re}z \leq a_1-r, \operatorname{Im}z=0\}$, $\gamma$ is contained in the only branch  $\{w=0\}$ of $E_1$ for $\{|z|>\frac{2}{3}\}$. Since both branching points of ${\Sigma}_1$ are contained in ${\Delta}_1=\{|z-a_1|<r\}$, and since only one of them is contained in the half-disk $\{z \in {\Delta}_1 :  \operatorname{Im}z>0\}$, we conclude that over the segment $\{a_1-r \leq \operatorname{Re}z \leq a_1+r, \operatorname{Im}z=0\}$ the curve $\gamma$ will  ``change from the branch $E_1$ to the branch ${\mathcal L}^{'}$''. Then, over the half-circle $\{|z-a_1|=r,  \operatorname{Im}z>0\}$ and the segment $\{a_2+r \leq \operatorname{Re}z \leq a_1-r, \operatorname{Im}z=0\}$ $ $, $\gamma$ is contained in  ${\mathcal L}^{'}$. After that, applying the same argument as we used for the segment  $\{a_1-r \leq \operatorname{Re}z \leq a_1+r, \operatorname{Im}z=0\}$, we conclude that, over the segment $\{a_2-r \leq \operatorname{Re}z \leq a_2+r, \operatorname{Im}z=0\}$, the curve $\gamma$ will  ``change from the branch ${\mathcal L}^{'}$ to the branch ${\mathcal L}^{''}$''. Then, over the segment $\{a_3+r \leq \operatorname{Re}z \leq a_2-r, \operatorname{Im}z=0\}$ and the half-circle $\{|z-a_3|=r,  \operatorname{Im}z>0\}$,  $\gamma$ is contained in ${\mathcal L}^{''}$.  After that, the same argument as above shows that, over the segment $\{a_3-r \leq \operatorname{Re}z \leq a_3+r, \operatorname{Im}z=0\}$, the curve $\gamma$ will  ``change from the branch ${\mathcal L}^{''}$ to the branch $E_1$''.  And finally, over the segment $\{a_3+r \leq \operatorname{Re}z \leq \frac{2}{3}, \operatorname{Im}z=0\}$, the curve $\gamma$ is contained in the ``upper'' branch of $E_1$  up to the endpoint $(\frac{2}{3}, 0)$, where we meet the first part of the curve ${\gamma}$ which is (for $|z|<\frac{2}{3}$) contained in the ``lower'' branch of $E_1$.

For each $z_0 \in {\pi}_z(\gamma)$ and each $s>0$, consider the sets 
$${\Gamma}_s(z_0)=\{(z_0,w) : \min_{(z_0, w^{'}) \in \gamma}|w-w^{'}|=s\}$$
and
$${\Omega}_s(z_0)=\{(z_0,w) : \min_{(z_0, w^{'}) \in \gamma}|w-w^{'}|<s\}.$$

Then, for $s$ small enough, each set ${\Omega}_s(z_0)$ is the union of
finitely many (at most three) disks in $\{z_0\} \times {\mathbb C}_w$,
which are disjoint if $z_0$ is far enough from the circle
$\{|z|=\frac{2}{3}\}$, and is the union of two connected components,
one of which is a disk and the other one is the union of two disks
having nonempty intersection, if $\{|z_0|<\frac{2}{3}\}$ and $z_0$ is
close enough to the circle $\{|z|=\frac{2}{3}\}$. As $|z_0|
\rightarrow \frac{2}{3}$ from the side $\{|z|<\frac{2}{3}\}$, the
centers of the two disks constituting the second connected component of ${\Omega}_s(z_0)$ become closer to each other, and for $\{|z_0| \geq \frac{2}{3}\}$ this component becomes just one disk. Each set ${\Omega}_s(z_0)$ has  a natural orientation induced from ${\mathbb C}_w$ and, hence, ${\Gamma}_s(z_0)=b{{\Omega}_s(z_0)}$ has also a natural orientation.

Consider the set
$$T_s=\bigcup_{{z_0} \in {\pi}_z(\gamma)}{\Gamma}_s(z_0).$$

Since the curve $\gamma$ is piece-wise smooth, it follows from the
definition of ${\Gamma}_s(z_0)$ that the set $T_s$ is a piece-wise
smooth surface of dimension two in $\Delta \times {\mathbb C}_w$ with
the boundary on the above chosen circle ${\mathcal C}_s$. Moreover, since $\gamma$ is oriented, and since each set ${\Gamma}_s(z_0)$ is oriented, we can also orient the surface $T_s$. Topologically, $T_s$ is a torus with a disk removed, ${\mathcal C}_s$ being the boundary of this disk. Since the curve $\gamma \subset E$ is transversal to the $w$-direction, we conclude that $T_s \subset \Delta\times {\mathbb C}_w\setminus E$ for $s$ sufficiently small. This implies that the homology class $\left[{\mathcal C}_s\right]$ of the circle ${\mathcal C}_s$ in $H_1(\Delta\times {\mathbb C}_w\setminus E, {\mathbb  Z})$ is trivial.

Now we observe that, for each point $(z,w) \in E^{\operatorname{reg}}$, the circle ${\mathcal C}_s(z,w)=\{(z,w{'}) : |w-w{'}|=s\}$ is homological to zero, if $s>0$ is small enough. Indeed, since the set $E^{\operatorname{reg}}$ is connected, there is a smooth curve $\tilde{\gamma} \subset E^{\operatorname{reg}}$ connecting the points $A$ and $(z,w)$. Then, for $s>0$ small enough, the set ${\mathcal M}_s=\{(z,w^{'}) : |w-w{'}|=s, (z,w) \in \tilde{\gamma}\}$ is a smooth ``cylinder'' of dimension two which is contained in $\Delta\times {\mathbb C}_w\setminus E$ and has its boundary on  ${\mathcal C}_s(z,w)$ and ${\mathcal C}_s$. Therefore, the circles ${\mathcal C}_s(z,w)$ and ${\mathcal C}_s$ represent the same homology class in  $H_1(\Delta\times {\mathbb C}_w\setminus E, {\mathbb  Z})$. Since ${\mathcal C}_s$ is already proved to be homological to zero in $\Delta\times {\mathbb C}_w\setminus E$, it follows that ${\mathcal C}_s(z,w)$ is also homological to zero in  $\Delta\times {\mathbb C}_w\setminus E$.

Finally, let ${\mathcal C}$ be any smooth closed curve in
$\Delta\times {\mathbb C}_w\setminus E$. Then, there is a
2-dimensional disk ${\mathcal D}$ smoothly embedded into $\Delta\times
{\mathbb C}_w$ such that ${\mathcal C}=b{\mathcal D}$. We can assume
that the disk ${\mathcal D}$ is in general position, in particular,
that ${\mathcal D}$ intersects $E$ in finitely many points
$\{(z_p,w_p)\}^k_{p=1}$ which are contained in
$E^{\operatorname{reg}}$. Without loss of generality, we can also
assume that ${\mathcal D}$ is parallel to the $w$-direction in a
neighbourhood of each point $(z_p,w_p)$. Then the disks ${\mathcal
  D}_s(z_p,w_p)=\{(z_p, w^{'}) : |w_p-w{'}| \leq s\}$ are contained in
${\mathcal D}$ for $s>0$ small enough. Therefore, ${\mathcal
  C}=b{\mathcal D}$ is homological  to $\bigcup^k_{p=1}{b{\mathcal
    D}_s(z_p,w_p)}$ in $\Delta\times {\mathbb C}_w\setminus E$, the
homology being  ${\mathcal D} \setminus \bigcup^k_{p=1}{{\mathcal D}_s(z_p,w_p)}$. Since each circle ${\mathcal C}_s(z,w)=b{\mathcal D}_s(z_p,w_p)$ is already proved to be homological to zero in $\Delta\times {\mathbb C}_w\setminus E$, we conclude that ${\mathcal C}$ is also homological to zero. The proof of the claim is now completed.
\bx\\

\vspace{-0.1truecm}

As an application of Claim 1 we show the following property of the set $E$.

\vspace{-0.1truecm}

\paragraph{Claim 2.} {\it There exists a neighbourhood $U(E)$ of the set $E$ which does not contain any subset of $\bar{\Delta} \times {\mathbb C}_w$ defined by a Weierstrass pseudopolynomial.}\\

\vspace{-0.6truecm}

\paragraph{Proof.} Assume, to get a contradiction, that every
neighbourhood $U(E)$ of $E$ contains a subset defined by a Weierstrass
pseudopolynomial. For $R$ big enough consider the circle ${\mathcal
  C}_R=\left\{ z=0, |w| = R \right\} \subset {\Delta \times {\mathbb
    C}_w} \setminus E$ oriented counterclock-wise in the
$w$-variable. Then, in view of Claim 1, there is a 2-chain $S$ such
that $bS={\mathcal C}_R$ and $\operatorname{supp}S \subset {\Delta
  \times {\mathbb C}_w} \setminus E$. The last inclusion implies that
there exists a neighbourhood $U(E)$ of $E$ such that
$\operatorname{supp}S \cap U(E)=\emptyset$. By our assumption, there
is a subset $\tilde{E}$ of $U(E)$ which is defined by a Weierstrass
pseudopolynomial, i.e. it has the form (1) with $a_1(z), a_2(z),..., a_m(z)$ being continuous functions. Since  $\operatorname{supp}S \cap \tilde{E}=\emptyset$, the homology class $\left[ {\mathcal C}_R \right]$ of the circle ${\mathcal C}_R$ in $H_1(\Delta\times {\mathbb C}_w\setminus \tilde{E}, {\mathbb  Z})$ is trivial. Consider the continuous map $\Phi: \Delta\times {\mathbb C}_w \setminus \tilde{E} \rightarrow S^1$ defined by
\begin{equation}\label{20}
\Phi(z,w) = \frac{w^m+a_1(z)w^{m-1}+...+a_m(z)}{|w^m+a_1(z)w^{m-1}+...+a_m(z)|}.
\end{equation}
Then, on one hand, $\left[{\mathcal C}_R\right]=0$ in $H_1(\Delta\times {\mathbb C}_w \setminus \tilde{E}, {\mathbb  Z})$ and, hence, $\Phi_{\ast}\left(\left[{\mathcal C}_R\right]\right)=0$ in $H_1(S^1, {\mathbb  Z})$. On the other hand, the term $w^m$ in the numerator of formula (11) will dominate  for $(z,w) \in {\mathcal C}_R$, if $R$ is big enough. Therefore, the degree of the restriction of $\Phi$ to ${\mathcal C}_R$ (it is a map from $S^1$ to $S^1$) is equal to $m$. Hence, $\Phi_{\ast}\left(\left[{\mathcal C}_R\right]\right)=m\left[S^1\right] \neq 0$ in $H_1(S^1, {\mathbb  Z})$. This gives the desired contradiction and proves the claim.\bx\\

\bigskip
\parindent 0pt

\end{document}